\documentclass[11pt]{article}
\usepackage[utf8]{inputenc}
\usepackage{amsfonts,amssymb,amsmath,amsthm}
\usepackage{enumitem}

\usepackage{color} 
   \definecolor{cites}{rgb}{0.50 , 0.00 , 0.00}  
   \definecolor{urls} {rgb}{0.00 , 0.00 , 0.50}  
   \definecolor{links}{rgb}{0.00 , 0.00 , 0.50}   
\usepackage[
      colorlinks=true,   
      citecolor=cites,   
      urlcolor=urls,     
      linkcolor=links,   
      pdftitle={Half-line compressions and finite sections of Schrödinger operators with integer-valued potentials},  
      pdfauthor={Marko Lindner, Riko Ukena}, 
      pdfpagemode=UseOutlines, 
      pdfstartview=FitH,       
      bookmarksopen=false      
   ]{hyperref}

\newcommand\eps\varepsilon
\newcommand\ph\varphi
\newcommand\spec{{\rm spec}\,}  

\newcommand\spess{{\rm spec}_{\rm ess}\,}

\newcommand\Lim{{\rm Lim}}

\newcommand\im{{\rm im}}

\newcommand\flip[1]{{#1}^{\sim}}

\newcommand\C{{\mathbb C}}
\newcommand\R{{\mathbb R}}

\newcommand\Z{{\mathbb Z}}
\newcommand\N{{\mathbb N}}
\newcommand\I{{\mathbb I}}

\newtheorem{theorem}{Theorem}[section]
\newtheorem{lemma}[theorem]{Lemma}
\newtheorem{corollary}[theorem]{Corollary}
\newtheorem{proposition}[theorem]{Proposition}
\newtheorem{definition}[theorem]{Definition}

\newenvironment{example}
 {\par\noindent\refstepcounter{theorem}{\bf Example \thetheorem}\ }
 {\raisebox{1mm}{\framebox{}}\pagebreak[2]}

\newenvironment{Proof}[1]{%
  \begin{proof}[{\bf Proof{#1}}]%
}{
   \end{proof}
}

\numberwithin{figure}{section}  

\newcounter{abccounter}

\makeatletter
\let\@fnsymbol\@arabic
\makeatother

\begin{document}
\title{\bf Half-line compressions and finite sections of discrete Schrödinger operators with integer-valued potentials}
\author{
{\sc Marko Lindner} \quad and \quad
{\sc Riko Ukena}\footnote{TU Hamburg (TUHH), Maths Institute, D-21073 Hamburg, Germany,  {\tt \{lindner,riko.ukena\}@tuhh.de}}
}

\date{\today}
\maketitle
\begin{quote}
\renewcommand{\baselinestretch}{1.0}
\footnotesize {\sc Abstract.}
We study 1D discrete Schrödinger operators $H$ with integer-valued potential and show that, $(i)$, invertibility (in fact, even just Fredholmness) of $H$ always implies invertibility of its half-line compression $H_+$ (zero Dirichlet boundary condition, i.e.~matrix truncation). In particular, the Dirichlet eigenvalues avoid zero -- and all other integers.
We use this result to conclude that, $(ii)$, the finite section method (approximate inversion via finite and growing matrix truncations) is applicable to $H$ as soon as $H$ is invertible.
The same holds for $H_+$.
\end{quote}

\noindent
{\it Mathematics subject classification (2020):} 47N40; Secondary 47B36, 47B93, 65J10.\\
{\it Keywords and phrases:} Schrödinger operator, half-line, finite section, spectrum 
\section{Introduction}
{\bf Discrete Schrödinger operators.\ }
We look at so-called {\sl discrete Schrödinger operators} in 1D, each acting via
\begin{equation} \label{eq:H}
(Hx)_n\ =\ x_{n-1} + v(n)x_n+ x_{n+1},\quad n\in\Z,
\end{equation}
as a bounded linear operator $H$ on $\ell^2(\Z)$. The (bounded) function $v:\Z\to\R$ is referred to as {\sl the potential} of $H$. The operator $H$ acts via matrix-vector multiplication by a two-sided infinite matrix $(H_{ij})_{i,j\in\Z}$ with main diagonal $H_{ii}=v(i)$, super and subdiagonal $H_{i,i\pm 1}=1$ and all other entries equal to zero.

Now let $\Z_+:=\{k\in\Z:k\ge 0\}$ and, given the operator \eqref{eq:H} on $\ell^2(\Z)$, we refer to the corresponding operator
\begin{equation} \label{eq:H+}
(H_+x)_n\ =\ x_{n-1} + v(n)x_n+ x_{n+1},\quad n\in\Z_+,\qquad\text{where}\quad x_{-1}:=0,
\end{equation}
on $\ell^2(\Z_+)$ as the {\sl half-line compression} $H_+$ of $H$.

From the operator perspective, $H_+$ equals $PHP|_{\im\, P}:\ell^2(\Z_+)\to\ell^2(\Z_+)$, where $P$ is the operator of multiplication by the characteristic function of $\Z_+$ (i.e.~the orthogonal projection of $\ell^2(\Z)$ onto $\ell^2(\Z_+)$). 
From the matrix perspective, $H_+$ corresponds to the one-sided infinite submatrix $(H_{ij})_{i,j\in\Z_+}$ of the two-sided infinite matrix $(H_{ij})_{i,j\in\Z}$ behind $H$.

In the same style, we denote the compression of $H$ to $\ell^2(\Z_-)$ by $H_-$, where $\Z_-:=-\Z_+$, the Dirichlet condition is $x_1=0$, and the corresponding matrix is $(H_{ij})_{i,j\in\Z_-}$.
\medskip

{\bf Finite sections.\ }
For $l,r\in\Z$ with $l<r$, let $H_{l..r}$ denote the compression of $H$ to the (space of functions on the) interval $\{l,\dots,r\}$ -- again with homogeneous Dirichlet conditions. The matrix behind $H_{l..r}$ is $(H_{ij})_{i,j=l}^r$, which is a so-called {\sl finite section} of $H$. The {\sl finite section method (FSM)} consists in approximating $H$ by a sequence $(H_n)$ of finite sections
\begin{equation} \label{eq:Hn}
H_n\ :=\ H_{l_n..r_n}\qquad\text{with}\qquad l_n\to -\infty\quad\text{and}\quad r_n\to +\infty
\end{equation}
for the asymptotic inversion (i.e.~the approximation of $H^{-1}$) or the spectral approximation of $H$. For a fixed choice of cut-off sequences $(l_n)$ and $(r_n)$, we will simply call the sequence $(H_n)$ from \eqref{eq:Hn} itself the FSM of $H$.
For an operator $H$ on $\ell^2(\Z_+)$ the FSM with cut-off sequence $(r_n)$, where $r_n\to+\infty$, is defined via $H_n\ :=\ H_{0..r_n}$. 

The FSM $(H_n)$ is called {\sl applicable to $H$} if all but finitely many $H_n$ are invertible and $H_n^{-1}b\to H^{-1}b$ for all $b\in\ell^2(\Z)$. Here, $H_n^{-1}$ is extended by zero, and the invertibility of $H$ is implicitly assumed as well. In particular, the FSM, if applicable, can be used to approximately solve equations $Hx=b$.

By a standard result of numerical analysis (``consistency+stability=convergence'' a.k.a.~Lax equivalence theorem~\cite{Lax} a.k.a.~Polski's theorem \cite{Polski}, see e.g.~\cite{BoeSi2,HaRoSi2}), one has
\begin{equation} \label{eq:FSMappl}
\left.
\begin{array}{r}
\text{the FSM $(H_n)$}\\
\text{is applicable to $H$}
\end{array}
\right\} \quad \iff\quad \left\{
\begin{array}{l}
\text{(a)\ \ $H$ is invertible,}\\
\text{(b)\ \ all but finitely many $H_n$ are invertible, and}\\
\text{(c)\ \ all inverses are uniformly bounded.}
\end{array}\right.
\end{equation}
Practically, condition (a) is mandatory anyway as assumption for the unique solvability of the equation $Hx=b$, and conditions (b) and (c) are typically a bit annoying to check.

The study of banded Toeplitz operators \cite{GohbergFeldman,BoeSi2} was a first instance where an operator class (on the half-line, finitely many constant diagonals, not limited to the tridiagonal setting \eqref{eq:H+} with $v\equiv const$) was identified for which condition (a) always implies (b) and (c), so that the applicability of the FSM already follows from the invertibility of the operator. Of course, such a situation simplifies FSM matters a lot: let us call $H$ {\sl FSM-simple} if condition (a) implies (b) and (c). For a particular operator $H$, this means, the FSM is applicable or $H$ is not even invertible.
\medskip

{\bf Example: The Fibonacci Hamiltonian.\ }
The standard 1D model for the study of electrical conduction properties of quasicrystals is the so-called {\sl Fibonacci Hamiltonian}. It is the discrete Schrödinger operator \eqref{eq:H} with potential
\begin{equation} \label{eq:FibPot}
v(n)\ =\ \chi_{[1-\alpha,1)}(n\alpha \!\!\mod  1)\, ,\qquad n\in\Z\,,
\end{equation}
where $\chi_I$ refers to the characteristic function of the interval $I$ and $\alpha=\frac12(\sqrt5-1)$ is the golden ratio.
Since $\alpha$ is irrational, $v$ is not periodic -- but it has very interesting combinatorial features instead. In~\cite{LiSoeding}, using many of those properties, it is shown that the FSM is applicable to the Fibonacci Hamiltonian, whence it is FSM-simple.
\medskip

{\bf Our results.\ }
Already then, Albrecht Böttcher sparked the discussion whether all those arguments using the complicated structure of \eqref{eq:FibPot} were really necessary in~\cite{LiSoeding} and what we would say if one day it turned out that every discrete Schrödinger operator \eqref{eq:H} with a $\{0,1\}$-valued (or even just integer) potential is FSM-simple.

We tried to give periodic counter-examples right away, failed to produce any and recently proved in \cite{Weber,perSchr} that none exist. In the current paper, as our first result, we even prove that periodicity has nothing to do with this non-existence; in short, Böttcher was right: 
\begin{theorem} \label{thm:FSMsimple}
Let $H$ be a discrete Schrödinger operator \eqref{eq:H} on $\ell^2(\Z)$ with an integer-valued potential $v:\Z\to\Z$. Then both $H$ and $H_+$ are FSM-simple.
\end{theorem}
Of course, this makes some of the twisted Fibonacci arguments in~\cite{LiSoeding} redundant although, for a FSM-simple operator, one is still left with proving invertibility in order to conclude applicability of the FSM.

Conditions (b) and (c) in \eqref{eq:FSMappl}, and hence the applicability of the FSM to $H$,
are connected to the invertibility of half-line compressions of so-called limit operators of $H$, see Lemma \ref{lem:FSMlimops} below. That's why the proof of Theorem \ref{thm:FSMsimple} relies on our second main result which is, moreover, interesting in its own right:
\begin{theorem} \label{thm:H+}
For a discrete Schrödinger operator $H$ on $\ell^2(\Z)$, as in \eqref{eq:H}, with an integer-valued potential $v:\Z\to\Z$, and its half-line compression $H_+$, as in \eqref{eq:H+}, the following implication holds:
\begin{equation} \label{eq:impl}
H\text{ is a Fredholm operator on }\ell^2(\Z)\qquad\implies\qquad H_+\text{ is invertible on }\ell^2(\Z_+).
\end{equation}
In particular, $H_+$ is invertible if $H$ is invertible.
\end{theorem}

{\bf Essential spectra and Dirichlet eigenvalues. }
In many situations, including periodic, almost-periodic, Sturmian and pseudo-ergodic potentials, the arrow in implication \eqref{eq:impl} points in the other direction (even stronger, with $H_+$ Fredholm $\implies$ $H$ invertible). More precisely, this is the case when $H$ is a limit operator of $H_+$. In short, this means that every finite subword of the potential $v$ of $H$ (understood as an infinite string) occurs infinitely often (at least up to arbitrary precision) in the right half of $v$.

Under this condition, by Lemma \ref{lem:spessH+} below, the spectrum of $H$ is equal to both the essential spectra of $H$ and $H_+$,
\begin{equation} \label{eq:spessH+}
\spec H\ =\ \spess H\ =\ \spess H_+\,.
\end{equation}
The spectrum of $H_+$ is typically a larger set. It arises from \eqref{eq:spessH+} by adding eigenvalues of $H_+$, the so-called {\sl Dirichlet-eigenvalues}. Their name addresses their cause: the truncation, i.e.~introduction of a homogeneous Dirichlet boundary condition.

If one is actually trying to compute $\spec H$ via a truncation technique like the FSM then such Dirichlet eigenvalues typically appear and lead to so-called {\sl spectral pollution}, being, erroneously, caused by the method rather than by the physical problem behind the operator $H$.

For example, for a $p$-periodic potential, one can show (see \cite{perSchr} and the references there) that \eqref{eq:spessH+} consists of at most $p$ closed intervals. Dirichlet eigenvalues of $H_+$ can only occur in the closure of the gaps in \eqref{eq:spessH+}, and each gap contains at most one Dirichlet eigenvalue. Similar results (see \cite{Dirichlet}) extend to aperiodic limits like the Fibonacci Hamiltonian, where \eqref{eq:spessH+} is a Cantor set, so that the number of gaps is infinite.

Remarkably, by our Theorem \ref{thm:H+}, zero is never among the Dirichlet eigenvalues, as well as all other integers $z\in\Z$, which is seen by applying the same result to $H-zI$.
    
{\bf Outline of the paper. }
After introducing all the necessary tools and language in Section \ref{sec:tools}, we prove both our theorems in Section \ref{sec:proofs} of the paper. The final section, Section \ref{sec:extend} shows some possible directions of extension but also some examples showing the limits of what is possible. 

\pagebreak

\section{Notations and tools} \label{sec:tools}
{\bf Spaces and operators. }
Let $\I\in\{\Z,\Z_+,\Z_-\}$ here and in what follows.
We demonstrate our results in the simplest case, $\ell^2(\I)$,
the standard space of all complex-valued sequences $x=(x_k)_{k\in\I}$ for which $\|x\|^2=\sum_{k\in\I}|x_k|^2<\infty$. Extension to $\ell^p$ with $p\in [1,\infty]$ is possible, see Section \ref{sec:extend}.

The statements in this section are not limited to discrete Schrödinger operators \eqref{eq:H} but apply to all so-called {\sl band operators} on $\ell^2(\I)$: operators $A$ whose matrix representation $(A_{ij})_{i,j\in\I}$ has uniformly bounded entries and support on finitely many diagonals only. Such an operator $A$ is always bounded on $\ell^2(\I)$. All this is clearly the case for our discrete Schrödinger operators $H$, $H_+$ and $H_-$.

{\bf Fredholmness, essential spectrum and limit operators. }
Recall that a bounded linear operator $A$ on a Banach space $X$ is a {\sl Fredholm operator}
if its coset, $A+K(X)$, modulo compact operators $K(X)$, is invertible in the so-called Calkin algebra $L(X)/K(X)$. This holds if and only if the nullspace of $A$ has finite dimension and the range of $A$ has finite codimension in $X$. In particular, Fredholm operators have a closed range. 
The {\sl essential spectrum} of $A$ is the spectrum of $A+K(X)$ in the Calkin algebra $L(X)/K(X)$, i.e.~the set of all $\lambda\in\C$ for which $A-\lambda I$ is not a Fredholm operator.

Since the coset $A+K(X)$ cannot be affected by changing finitely many matrix entries, its study takes place ``at infinity''. This is where limit operators \cite{RaRoSi1998,RaRoSiBook,Li:Book} come in:

\begin{definition}
For a band operator $A$ on $\ell^2(\I)$, we look at all its translates $S^{-k}AS^k$ with $k\in\Z$ and speak of a {\sl limit operator}, $A_h$, on $\ell^2(\Z)$ if, for a particular sequence $h=(h_n)$ in $\Z$ with $|h_n|\to\infty$, the corresponding sequence of translates, $S^{-h_n}AS^{h_n}$, converges strongly to $A_h$.

Hereby, $(S x)_n = x_{n - 1}$ denotes the shift operator on $\ell^2(\Z)$ and a sequence $A_n$ is said to converge strongly to $A$ if $A_nx\to Ax$ 
for all $x\in\ell^2(\Z)$.
Moreover, let $\Lim(A)$ denote the set of all limit operators of $A$, together with the following local versions: For an integer sequence $g=(g_n)$ with $|g_n|\to \infty$, we put 
$$
\Lim_g(A)\ :=\ \{A_h\ :\ h\text{ is a subsequence of } g\}
$$
as well as $\Lim_+(A):=\Lim_{(1,2,\dots)}(A)$ and $\Lim_-(A):=\Lim_{(-1,-2,\dots)}(A)$.
\end{definition}

From the matrix perspective, 
$A_h=(\tilde A_{ij})_{i,j\in\Z}$ is the limit operator of $A=(A_{ij})_{i,j\in\I}$ with respect to the sequence $h=(h_n)$ in $\Z$ if, for all $i,j\in\Z$,
\begin{equation}\label{eq:limop}
A_{i+h_n,j+h_n}\to \tilde A_{ij}\quad\textrm{as}\quad n\to\infty.
\end{equation}
Here is the announced connection to Fredholm operators.

\begin{lemma} \label{lem:Fredholm}
For a band operator $A$ on $\ell^2(\I)$, it holds that the following are equivalent:\\[-7mm]
\begin{enumerate}[label=(\alph*)]\itemsep-1mm
\item $A$ is a Fredholm operator on $\ell^2(\I)$,
\item  all limit operators of $A$ are invertible on $\ell^2(\Z)$ \cite{RaRoSi1998,LiSei:BigQuest}.
\end{enumerate}
\end{lemma}
As a consequence of Lemma \ref{lem:Fredholm},
\begin{equation} \label{eq:specLim}
\spess A\ =\ \bigcup_{A_h\in\Lim(A)} \spec A_h\,.
\end{equation}
Let us generalize the definition of a half-line compression $A_+:=PAP:\ell^2(\Z_+)\to\ell^2(\Z_+)$, and similarly $A_-$, from below \eqref{eq:H+} to band operators $A$ on $\ell^2(\Z)$. Then we have the following result leading to formula \eqref{eq:spessH+} above.
\begin{lemma} \label{lem:spessH+}
If a band operator $A$ on $\ell^2(\Z)$ is limit operator of its half-line compression $A_+$, in short, $A\in\Lim(A_+)$, then
\[
\spec A\ =\ \spess A\ =\ \spess A_+\,.
\]
\end{lemma}
\begin{Proof}{}
The inclusion $\spess A_+\subset \spess A$ holds by \eqref{eq:specLim} since $\Lim(A_+)=\Lim_+(A)\subset\Lim(A)$.
The inclusion $\spess A\subset\spec A$ is standard, of course.
Finally, $\spec A\subset \spess A_+$ holds by \eqref{eq:specLim} and $A\in\Lim(A_+)$, see our assumption.
\end{Proof}

{\bf Limit operators and the FSM. }
%
%
%
%
%
Fix integer sequences $l=(l_n)_{n\in\N}$ and $r=(r_n)_{n\in\N}$ with $l_n < r_n$ for all $n\in\N$ and $l_n \to -\infty$ and $r_n \to +\infty$ as $n\to\infty$. For a band operator $A$ on $\ell^2(\Z)$, resp.~$\ell^2(\Z_+)$, look at its sequence $(A_n)_{n\in\N}$ of finite sections
\[
A_n\ :=\ A_{l_n..r_n}\ =\ (A_{ij})_{i,j=l_n}^{r_n},
\qquad\text{resp.}\qquad
A_n\ :=\ A_{0..r_n}\ =\ (A_{ij})_{i,j=0}^{r_n},
\qquad n\in\N.
\]
Recall from our introduction that we call the sequence $(A_n)$ itself the finite section method (FSM) of $A$ and say that it is applicable if the conditions (a), (b), (c) hold in \eqref{eq:FSMappl}.
\begin{lemma}{\cite{RaRoSi:FSMsubs,Li:FSMsubs,LindnerRoch}} \label{lem:FSMlimops}
{\bf a) } The two-sided case:
The FSM with cut-off sequences  $l=(l_n)_{n\in\N}$ and $r=(r_n)_{n\in\N}$ as above is applicable to a band operator $A$ on $\ell^2(\Z)$ if and only if the following operators are invertible:\\[-7mm]
 \begin{enumerate}[label=(\alph*)]\itemsep-1mm
  \item the operator $A$ itself,
  \item all operators $L_+$ with $L\in\Lim_l(A)$,
  \item all operators $R_-$ with $R\in\Lim_r(A)$.
 \end{enumerate}
 
 {\bf b) } The one-sided case:
 The FSM with cut-off sequence $r=(r_n)_{n\in\N}$ as above is applicable to a band operator $A_+$ on $\ell^2(\Z_+)$ if and only if the following operators are invertible:\\[-7mm]
 \begin{enumerate}[label=(\alph*)]\itemsep-1mm 
 \setcounter{enumi}{3}
  \item the operator $A_+$ itself,
  \item all operators $R_-$ with $R\in\Lim_r(A_+)$.
 \end{enumerate}
\end{lemma}

Without loss, we can restrict ourselves to the study of half-line compressions in the positive direction. Indeed, by an elementary reflection technique, one can connect the study of $R_-$ with $R\in\Lim_r(A)$ to the study of 
\[
\flip R_+\ :=\ (\flip R)_+,\quad\text{ where }\quad \flip R:=\Phi R\Phi
\]
with the flip operator $(\Phi x)_n=x_{-n}$ for all $n\in\Z$.
In matrix language, $\flip R=(R_{-i,-j})_{i,j\in\Z}$ if $R=(R_{ij})_{i,j\in\Z}$. Note that $\flip R\in\Lim_{-r}(\flip A)$.
It is straightforward to check that $R$ and $\flip R$ are simultaneously invertible on $\ell^2(\Z)$ and that $R_-$ is invertible on $\ell^2(\Z_-)$ if and only if $\flip R_+$ is invertible on $\ell^2(\Z_+)$.
%
%


{\bf Discrete Schrödinger operators and transfer matrices. }
In order to examine the kernel of $H$ or $H_+$, it is useful to reformulate the scalar three-term recurrence
\begin{equation*}
  0\ =\ (Hx)_n \ =\ x_{n-1} + v(n)x_n + x_{n+1}, \quad n\in\Z 
\end{equation*}
as a vector-valued two-term recursion 
\begin{equation} \label{eq:vecrecur}
 \begin{pmatrix}
  x_n \\ x_{n+1}
 \end{pmatrix}
 = 
 \begin{pmatrix}
   0 & 1 \\ -1 & -v(n) 
 \end{pmatrix}
\begin{pmatrix}
 x_{n-1} \\ x_n
\end{pmatrix} \,.
\end{equation}
The matrix $T(n) := \begin{pmatrix}   0 & 1 \\ -1 & -v(n)  \end{pmatrix} $ from \eqref{eq:vecrecur} is the so-called (and well-known) {\sl transfer matrix}.

\section{The proofs of our results} \label{sec:proofs}
The proof of Theorem \ref{thm:FSMsimple} uses Theorem \ref{thm:H+}. The latter follows from the following two results:
\begin{proposition} \label{prop:int->inj}
Let $H_+$ be a discrete Schrödinger operator on the half-line, $\ell^2(\Z_+)$, as in \eqref{eq:H+}, with an integer-valued potential $v:\Z_+\to\Z$. Then $H_+$ is injective.
\end{proposition}
\begin{Proof}{}
In order to show that $H_+$ is injective, take a vector $x\in\ell^2(\Z_+)$ with $H_+ x = 0$.
We further assume that $x\neq 0$, which will lead to a contradiction.
Recall from \eqref{eq:H+} that the Dirichlet boundary condition $x_{-1}= 0$ holds.
Since $x_0=0$ would imply that $x=0$, we can, without loss of generality, assume that $x_0=1$.
By \eqref{eq:vecrecur} and the definition of the transfer matrices $T(k)$, the other entries of $x$ satisfy the condition
\begin{equation}\label{eq:transfers2}
\begin{pmatrix}
  x_n \\ x_{n+1}
 \end{pmatrix}
 \ =\  
 T(n) \dots T(0) 
\begin{pmatrix}
 x_{-1} \\ x_0
\end{pmatrix}
\ =\  
 T(n) \dots T(0) 
\begin{pmatrix}
 0 \\ 1
\end{pmatrix},\qquad n\in\N.
\end{equation}
Since, by assumption, $v(k)\in\Z$ for all $k\in\Z_+$, we find that  $T(k) = \begin{pmatrix}   0 & 1 \\ -1 & -v(k)  \end{pmatrix} \in \Z^{2\times 2}$.
Thus \eqref{eq:transfers2} implies that $x_n\in\Z$ for all $n\in\Z_+$. 
As an integer sequence in $\ell^2(\Z_+)$, $x$ can only have finitely many non-zero entries.
Thus there exists an $n\in\N$ such that $x_n=x_{n+1}=0$.
Using the invertibility of the transfer matrices, we get
\begin{equation}\label{eq:transfers3}
\begin{pmatrix}
  x_{-1} \\ x_{0}
 \end{pmatrix}
 = 
 T(0)^{-1} \dots T(n)^{-1} 
\begin{pmatrix}
 x_n \\ x_{n+1}
\end{pmatrix} =
\begin{pmatrix}
 0 \\ 0
\end{pmatrix}
\end{equation}
and hence, $x_0=0$. This contradicts the assumption $x_0\neq 0$.
Hence we have shown that $H_+$ cannot have any nontrivial vectors $x$ with $H_+x=0$, i.e.~it is injective.
\end{Proof}
So, after this proposition, all that is missing for the invertibility of the self-adjoint operator $H_+:\ell^2(\Z_+)\to\ell^2(\Z_+)$ with integer-valued potential is the closedness of its range.
\begin{lemma} \label{lem:closedrange}
Let $H$ be a discrete Schrödinger operator, \eqref{eq:H}, on $\ell^2(\Z)$ and $H_+$ its half-line compression, \eqref{eq:H+}, on $\ell^2(\Z_+)$. Then the following implications hold,
\[
(i)\iff(ii)\implies (iii)\iff(iv),
\]
where\\
\begin{tabular}{p{15mm}rlp{5mm}rl}
&$(i)$ & $H$ is a Fredholm operator, & &
$(iii)$ & $H_+$ is a Fredholm operator,\\
&$(ii)$ & $H$ has a closed range, & &
$(iv)$ & $H_+$ has a closed range.
\end{tabular}
\end{lemma}
\begin{Proof}{}
$(i)\Rightarrow (ii)$ and $(iii)\Rightarrow (iv)$ are standard.
Their reverse directions follow, once we show that always $\dim\ker(H)\le 2$ and $\dim\ker(H_+)\le 1$.
For $H_+$ this is shown in the proof of Proposition \ref{prop:int->inj}: A vector $x$ in the kernel
of $H_+$ is given by $x_{-1}=0$, $x_0=a$ and all other $x_n$ by \eqref{eq:transfers2}. So we have one degree of freedom, $a$, (or possibly just $\ker(H_+)=\{0\}$). For $x\in\ker(H)$, one has two degrees of freedom, say $x_{-1}$ and $x_0$, and the recurrence \eqref{eq:transfers2} is going forward and (appropriately inverted) backward from there.

The implication $(i)\Rightarrow (iii)$ holds by Lemma~\ref{lem:Fredholm} and $\Lim(H)\supset\Lim_+(H)= \Lim(H_+)$.
\end{Proof}

\begin{Proof}{\ of Theorem \ref{thm:H+}}
Let $H$ be a discrete Schr\"odinger operator on $\ell^2(\Z)$ with an integer-valued potential $v:\Z\to\Z$ and let $H$ be a Fredholm operator.
Since $H_+$ is self-adjoint, it suffices to show that it is a) injective and has b) a closed range.\\[-7mm]
\begin{itemize} \itemsep-1mm
\item[a)] follows by Proposition \ref{prop:int->inj} and $v(n)\in\Z$ for all $n\in\Z$,
\item[b)] follows from Lemma \ref{lem:closedrange}.
\qedhere
\end{itemize}
\end{Proof}

\begin{Proof}{\ of Theorem \ref{thm:FSMsimple}}
We start with the two-sided case, $H$ on $\ell^2(\Z)$.
Let $H$ be invertible. By Lemma~\ref{lem:Fredholm}, all $L\in \Lim_-(H)$ and $R\in \Lim_+(H)$ are also invertible. Because the convergence in \eqref{eq:limop} is by being eventually constant, $L$ and $R$ are again discrete Schrödinger operators with integer potential.

Clearly, also $\flip R$ is invertible and a discrete Schrödinger operator with integer potential.

Applying Theorem \ref{thm:H+} to $H=L$ and $H=\flip R$, we get that $L_+$ and $\flip R_+$, and consequently also $R_-$, is invertible on the corresponding half-line space. 

Because this is true for all $L\in \Lim_-(H)$ and $R\in \Lim_+(H)$,
we get, by Lemma \ref{lem:FSMlimops} a), that the FSM with arbitrary cut-off sequences
$l=(l_n)_{n\in\N}$ and $r=(r_n)_{n\in\N}$ with $l_n < r_n$ for all $n\in\N$ and $l_n \to -\infty$ and $r_n \to +\infty$ as $n\to\infty$ is applicable to $H$.

In the one-sided case, $H_+$ on $\ell^2(\Z_+)$, assume that $H_+$ is invertible. Now argue with Lemma~\ref{lem:FSMlimops} b) exactly as above (without the operators $L_+$) to see that the FSM is applicable to $H_+$.
\end{Proof}

\begin{corollary} \label{cor:final}
For a discrete Schrödinger operator $H$ on $\ell^2(\Z)$, as in \eqref{eq:H}, with an integer-valued potential $v:\Z\to\Z$, and its half-line compression $H_+$, as in \eqref{eq:H+}, the following holds:

If $H$ has a closed range then the FSM is applicable to $H_+$.
\end{corollary}
\begin{Proof}{}
$\im(H)$ closed $\stackrel{\rm L \ref{lem:closedrange}}\implies$
$H$ Fredholm $\stackrel{\rm T \ref{thm:H+}}\implies$
$H_+$ invertible $\stackrel{\rm T \ref{thm:FSMsimple}}\implies$
FSM applicable to $H_+$.
\end{Proof}

%
\section{Limitations and possible directions of extension} \label{sec:extend}
{\bf Limitations. }\\[-6mm]
\begin{itemize}
\item We cannot replace ``integer-valued'' by ``rational-valued'' in Theorems \ref{thm:FSMsimple}, \ref{thm:H+} and Proposition \ref{prop:int->inj}, see the following example:
\end{itemize}
\begin{example}
Let $v$ be the $3$-periodic, two-sided extension of the vector $(\frac 12, 2,\frac 12)$ with $v(1)=2$. Then it is shown in \cite[Example 4.2]{perSchr} that the corresponding operator $H$ from \eqref{eq:H} is invertible, while $H_+$ from \eqref{eq:H+} is not invertible. $H_+$ is Fredholm, though: note that $\Lim(H_+)=\{H,\,S^{-1}HS,\,S^{-2}HS^2\}$. So $H_+$ must lack injectivity. 
In particular, the claims of Proposition \ref{prop:int->inj} and Theorem \ref{thm:H+}, do obviously not apply in this situation. Theorem \ref{thm:FSMsimple}: also not, see \cite{perSchr}.
\end{example}

\begin{itemize}
\item Unlike $H_+$, the full-space operator $H$ is not automatically injective if it has an integer-valued potential. See the following example:
\end{itemize}
\begin{example}
We construct a potential $v\in\{0,1\}^\Z$ with $v_\pm:=v|_{\Z_\pm}$ periodic each.
$v_+$ first: In the sense of \eqref{eq:transfers2}, we get from ${x_{-1}\choose x_0}$ to ${x_{kp-1}\choose x_{kp}}$ by applying $M^k$ with $M:=T(p-1)\dots T(0)$, where $k\in\N$ and $p$ is the period of $v_+$. Then $x_{kp}=\alpha \lambda_1^k+\beta\lambda_2^k$, where $\lambda_i\in\spec M$ and $\alpha,\beta$ depend on ${x_{-1}\choose x_0}$, e.g.~\cite{perSchr}.
Noting that $\lambda_1\lambda_2=\det M=\prod_{k=0}^{p-1} \det T(k)=1$, we choose $p=5$ and $v|_{\{0,\dots,4\}}=10101$, so that $|\lambda_1|<1$ and ${x_{-1}\choose x_0}\ne 0$ has no component in the eigenspace of $\lambda_2$, whence $\beta=0$ and $x_{kp}$ (and hence $x_n$) decays exponentially as $k,n\to+\infty$.
Now choose $v_-$ with period $q$ so that also 
\[
{x_{-kq-1}\choose x_{-kq}}=\Big(T(-q)^{-1}\dots T(-1)^{-1}\Big)^k{x_{-1}\choose x_0},\qquad k\in\N
\]
decays exponentially. The trick is to take $v_-$ as reflection of $v_+$ and replace every $0$ by $000$ and every $1$ by $11$. The result is $x_{-kq}=x_{kp}$ since $T_0=T_0^{-3}$ and $T_1=T_1^{-2}$ for $T_\gamma={0~~1\choose -1~-\gamma}$.
Consequently, $q=12$, $v|_{\{-12,\dots,-1\}}=110001100011$ and $Hx=0$ with $x\in\ell^2(\Z)\setminus\{0\}$.
\end{example}
\medskip 

\begin{itemize}
\item The majority of our results (all but Lemma \ref{lem:closedrange}) have only very restricted counterparts in the non-selfadjoint situation
\begin{equation} \label{eq:HNSA}
(Hx)_n\ =\ ax_{n-1} + v(n)x_n+ bx_{n+1},\quad n\in\Z,
\end{equation}
with $a,b\in\C$, where  $T(n) = \begin{pmatrix}   0 & 1 \\ -\frac ab & -\frac{v(n)}b  \end{pmatrix}$.

Precisely, all results that require an integer potential do in fact rely on Proposition \ref{prop:int->inj} and integer transfer matrices $T(n)\in\Z^{2\times 2}$, which then translates to $a$ and all $v(n)$ being integer multiples of $b$.

However, for both our theorems, we need, besides injectivity and a closed range of the half-line compression, also injectivity of its adjoint (which is now different). 
Arguing via Proposition~\ref{prop:int->inj} again,
in order to guarantee integer transfer matrices, also $\frac ba$ and all $\frac{v(n)}a$ have to be integer.
So in total, all $v(n)$ have to be integer multiples of $a=\pm b\ne 0$. 

\end{itemize}

\noindent
{\bf Possible extensions. }\\[-8mm]
\begin{itemize} \itemsep-1mm
\item The results are clearly not different when $H_-$ instead of $H_+$ is considered.

\item Also the starting point of the half-axis, here zero, is clearly irrelevant.

\item We have shown that also compressions to finite intervals (i.e.~finite sections), with homogeneous Dirichlet conditions on both ends, $l_n\to-\infty$ and $r_n\to +\infty$, inherit invertibility from $H$ for large enough $n$. 

Note that, in the general situation of all integer-valued potentials $v$, there is no interval length $L$ for which all finite sections $H_{l_n..r_n}$ are invertible as soon as they have sufficient size, $r_n-l_n+1=:L_n\ge L$. Indeed, for every odd $L\in\N$, $v$ can contain $L$ consecutive zeros, which makes a corresponding finite section of size $L_n=L$ singular (simple exercise by row operations), but $H$ can still be invertible.

For more specific classes, the situation is different: For the Fibonacci Hamiltonian $H$, it is shown in \cite{Weber} that $L=6$, and for its half-line compression $H_+$, $L=1$ (with $l_n\equiv 0$).

%
%

\item We can clearly pass from $\ell^2$ to $\ell^p$, even with $p\in[1,\infty]$.
Instead of the strong convergence of $S^{-h_n}AS^{h_n}$, one then looks at the so-called $\mathcal P$-convergence \cite{Sei:Survey}. 
\item We can generalize from potentials $v$ with values in $\Z$ to values in a set $R\subset\C$ with

\begin{tabular}{p{3mm}rlp{8mm}rl}
&$(i)$ & $-1, 0, 1\in R$, & &
$(iii)$ & $(R,+,\cdot)$ is a ring,\\
&$(ii)$ & $v(n)\in R$ for all $n\in \Z$, & &
$(iv)$ & $0$ is an isolated point of $R$.
\end{tabular}

From $(i), (ii), (iii)$ it follows that $T(n)\in R^{2\times 2}$ and hence $x_n\in R$ in \eqref{eq:transfers2} for all $n\in\Z_+$. From $(iv)$ it follows that, just as two lines below \eqref{eq:transfers2}, $x\in\ell^2(\Z_+)$ implies $x\in c_{00}(\Z_+)$.

By $(i)$ and $(iii)$, $s-t=s+(-1)\cdot t\in R$ for all $s,t\in R$. By $(iii)$, with every $r\in R$, also $r^2, r^3, \ldots\in R$, so that, by $(iv)$, the only $r\in R$ with $|r|<1$ is $r=0$. As a consequence of these observations, $R$ is discrete with $|s-t|\ge 1$ for all $s,t\in R$ with $s\ne t$.

Examples of such a set $R$ are grids $R=r^1\Z+r^2\Z+\dots+r^n\Z$ with $n\in\N$ and a fixed $r\in\C$ with $r^n=1$. In particular, this leads to $R=\Z$ for $n\in\{1,2\}$, $R=\Z+i\Z$ for $n=4$, while $n\in\{3,6\}$ lead to the same honeycomb grid $R$ (since $3$ is odd and $\Z=-\Z$). Note that $n>6$ is impossible as it leads to the existence of different $s,t\in R$ with $|s-t|<1$ and $n=5$ has the same problem as its grid coincides with the one for $n=10>6$.

\end{itemize}

\medskip

{\bf Acknowledgements.}
The authors thank Julian Großmann, Fabian Gabel, Dennis Gallaun, Daniel Lenz and Albrecht Böttcher for inspiring and helpful discussions.

\vfill
\noindent {\bf Authors' addresses:}\\
\\
Marko Lindner\hfill \href{mailto:lindner@tuhh.de}{\tt lindner@tuhh.de}\\
Riko Ukena\hfill \href{mailto:riko.ukena@tuhh.de}{\tt riko.ukena@tuhh.de}\\[2mm]
Institut Mathematik\\
TU Hamburg (TUHH)\\
D--21073 Hamburg\\
GERMANY

\end{document}